# ON HELICES AND BERTRAND CURVES IN EUCLIDEAN 3-SPACE

Murat Babaarslan[1] and Yusuf Yayli[2]

[1]Department of Mathematics, Faculty of Arts and Sciences
Bozok University, Yozgat, Turkey
murat.babaarslan@bozok.edu.tr
[2]Department of Mathematics, Faculty of Science
Ankara University, Tandogan, Ankara, Turkey
yayli@science.ankara.edu.tr

**Abstract**- In this article, we investigate Bertrand curves corresponding to the spherical images of the tangent, binormal, principal normal and Darboux indicatrices of a space curve in Euclidean 3-space. As a result, in case of a space curve is a general helix, we show that the curves corresponding to the spherical images of its the tangent indicatrix and binormal indicatrix are both Bertrand curves and circular helices. Similarly, in case of a space curve is a slant helix, we demonstrate that the curve corresponding to the spherical image of its the principal normal indicatrix is both a Bertrand curve and a circular helix.

**Key Words**- Helix, Bertrand curve, Sabban frame, Spherical images

## 1. INTRODUCTION

In the differential geometry of a regular curve in Euclidean 3-space, it is well known that, one of the important problems is characterization of a regular curve. The curvature $\kappa$ and the torsion $\tau$ of a regular curve play an important role to determine the shape and size of the curve [1, 2]. For example: if $\kappa = \tau = 0$, then the curve is a geodesic. If $\kappa \neq 0$ (constant) and $\tau = 0$, then the curve is a circle with radius $1/\kappa$.

A curve of constant slope or general helix in Euclidean 3-space is characterized by the property that the tangent lines make a constant angle with a fixed direction (the axis of the general helix). A classical result about helix stated by Lancret in 1802 and first proved by de Saint Venant in 1845 (see [3] for details) says that: A necessary and sufficient condition that a curve be a general helix is that the ratio $\kappa/\tau$ is constant along the curve, where $\tau \neq 0$. If both $\kappa$ and $\tau$ are non-zero constants, it is called a circular helix.

There are a lot of interesting applications of helices (e.g., DNA double and collagen triple helix, helical staircases, helical structures in fractal geometry and so on). All these make authors say that the helix is one of the most fascinated curves in science and nature.

In the study of fundamental theory and the characterizations of space curves, the corresponding relations between the curves are very interesting problem. The well-known Bertrand curve is characterized as a kind of such corresponding relation between the two curves. Bertrand curves discovered by J. Bertrand in 1850 are one of the important and interesting topics of classical special curve theory. A Bertrand curve is defined as a special curve which shares its the principal normals with another special curve, called Bertrand mate or Bertrand partner curve. The curve $\tilde{\gamma}$ is a Bertrand curve if and only if there exist non-zero real numbers $A$, $B$ such that $A\kappa(s) + B\tau(s) = 1$ for any $s \in I$ [1, 6]. So a circular helix is a Bertrand curve. Bertrand mates represent particular examples of offset curves [4] which are used in computer-aided design (CAD) and computer-aided manufacture (CAM).

Izumiya and Takeuchi [5] have shown that cylindrical helices can be constructed from plane curves and Bertrand curves can be constructed from spherical curves. After that, they [6] have studied cylindrical helices and Bertrand curves from the view point as curves on ruled surfaces. Schief [7] has given a study of the integrability of Bertrand curves. Kula and Yayli [8] have studied the spherical images of the tangent indicatrix and binormal indicatrix of a slant helix and they have shown that the spherical images are spherical helices. Camci et al. [9] have given some characterizations for a non-degenerate curve $\alpha$ to be a generalized helix by using harmonic curvatures of the curve in $n$-dimensional Euclidean space $\mathbb{E}^n$.

The aim of this paper is to investigate Bertrand curves corresponding to the spherical images of a space curve in Euclidean 3-space. In case of a space curve is general helix, we show that the curves corresponding to the spherical images of its the tangent indicatrix and binormal indicatrix are both Bertrand curves and circular helices. Furthermore, in case of a space curve is slant helix, we demonstrate that the curve corresponding to the spherical image of its the principal normal indicatrix is both a Bertrand curve and a circular helix.

## 2. PRELIMINARIES

To meet the requirements in the next section, here, the basic elements of the theory of curves in the Euclidean space are briefly presented. After that, we describe the method to construct Bertrand curves from spherical curves. A curve $\tilde{\gamma} : I \to \mathbb{R}^3$ with unit speed, is a general helix if there is some constant vector $\boldsymbol{u}$, so that $\boldsymbol{T} \cdot \boldsymbol{u} = \cos\theta$ is constant along the curve, where $\boldsymbol{T}(s) = \tilde{\gamma}'(s)$ and we call $\boldsymbol{T}(s)$ as a unit tangent vector of $\tilde{\gamma}$ at $s$. We define the curvature of $\tilde{\gamma}$ by $\kappa(s) = \|\tilde{\gamma}''(s)\|$. If $\kappa(s) \neq 0$, then the unit principal normal vector $\boldsymbol{N}(s)$ of the curve $\tilde{\gamma}$ at $s$ is given by $\tilde{\gamma}''(s) = \kappa(s)\boldsymbol{N}(s)$.
The unit vector $\boldsymbol{B}(s) = \boldsymbol{T}(s) \times \boldsymbol{N}(s)$ is called the unit binormal vector of $\tilde{\gamma}$ at $s$. For the derivations of the Frenet frame, the Frenet-Serret formulae hold:

$$\begin{bmatrix} \boldsymbol{T}'(s) \\ \boldsymbol{N}'(s) \\ \boldsymbol{B}'(s) \end{bmatrix} = \begin{bmatrix} 0 & \kappa(s) & 0 \\ -\kappa(s) & 0 & \tau(s) \\ 0 & -\tau(s) & 0 \end{bmatrix} \begin{bmatrix} \boldsymbol{T}(s) \\ \boldsymbol{N}(s) \\ \boldsymbol{B}(s) \end{bmatrix}, \qquad (1)$$

where $\tau(s)$ is the torsion of the curve $\tilde{\gamma}$ at $s$. For any unit speed curve $\tilde{\gamma} : I \to \mathbb{R}^3$, we call $\boldsymbol{W}(s) = \tau(s)\boldsymbol{T}(s) + \kappa(s)\boldsymbol{B}(s)$ as the Darboux vector field of $\tilde{\gamma}$. Spherical images (indicatrices) are well-known in classical differential geometry of curves. Let us define the curve $\boldsymbol{C}$ on $\mathbb{S}^2$ by the help of vector field $\boldsymbol{C}(s) = \boldsymbol{W}(s) / \|\boldsymbol{W}(s)\|$. This curve is called the spherical Darboux image or the Darboux indicatrix of $\tilde{\gamma}$. Similarly, the unit tangent vectors along the curve $\tilde{\gamma}$ generate a curve $(T)$ on $\mathbb{S}^2$. The curve $(T)$ is called the spherical indicatrix of $T$ or more commonly, $(T)$ is called the tangent indicatrix of the curve $\tilde{\gamma}$. If $\tilde{\gamma} = \tilde{\gamma}(s)$ is a natural representation of $\tilde{\gamma}$, then $(T) = T(s)$ will be a representation of $(T)$. One considers the principal normal indicatrix $(N) = N(s)$ and the binormal indicatrix $(B) = B(s)$ [3].

For a general parameter $t$ of a space curve $\tilde{\gamma}$, we can calculate the curvature and the torsion as follows:

$$\kappa(t) = \frac{\|\tilde{\gamma}'(t) \times \tilde{\gamma}''(t)\|}{\|\tilde{\gamma}'(t)\|^3}, \quad \tau(t) = \frac{\det(\tilde{\gamma}'(t), \tilde{\gamma}''(t), \tilde{\gamma}'''(t))}{\|\tilde{\gamma}'(t) \times \tilde{\gamma}''(t)\|^2}. \tag{2}$$

A curve $\tilde{\gamma}$ with $\kappa(s) \neq 0$ is called a slant helix if the principal normal lines of $\tilde{\gamma}$ make a constant angle with a fixed direction. $\tilde{\gamma}$ is a slant helix if and only if the geodesic curvature of the spherical image of the principal normal indicatrix $(N)$ of $\tilde{\gamma}$

$$\psi(s) = \left(\frac{\kappa^2}{(\kappa^2 + \tau^2)^{3/2}} \left(\frac{\tau}{\kappa}\right)'\right)(s) \tag{3}$$

is a constant function [10].

Now we give a new frame different from the Frenet frame. Let $\gamma : I \to \mathbb{S}^2$ be a unit speed spherical curve. We denote $\sigma$ as the arc-length parameter of $\gamma$. Let us denote $t(\sigma) = \gamma'(\sigma)$, and we call $t(\sigma)$ a unit tangent vector of $\gamma$ at $\sigma$. We now set a vector $s(\sigma) = \gamma(\sigma) \times t(\sigma)$, where $\gamma$ denotes the position vector of the curve. By definition of the curve $\gamma$, we have an orthonormal frame $\{\gamma(\sigma), t(\sigma), s(\sigma)\}$ along $\gamma$. This frame is called the Sabban frame of $\gamma$ [11]. Then we have the following spherical Frenet-Serret formulae of $\gamma$:

$$\begin{bmatrix} \gamma'(\sigma) \\ t'(\sigma) \\ s'(\sigma) \end{bmatrix} = \begin{bmatrix} 0 & 1 & 0 \\ -1 & 0 & \kappa_g(\sigma) \\ 0 & -\kappa_g(\sigma) & 0 \end{bmatrix} \begin{bmatrix} \gamma(\sigma) \\ t(\sigma) \\ s(\sigma) \end{bmatrix}, \tag{4}$$

where $\kappa_g(\sigma)$ is the geodesic curvature of the curve $\gamma$ on $\mathbb{S}^2$ which is given by $\kappa_g(\sigma) = \det(\gamma(\sigma), t(\sigma), t'(\sigma))$.

We now define a space curve

$$\tilde{\gamma}(\sigma) = a\int_{\sigma_0}^{\sigma} \gamma(\sigma)d\sigma + a\cot\theta \int_{\sigma_0}^{\sigma} s(\sigma)d\sigma + c, \qquad (5)$$

where $a$, $\theta$ are constant numbers and $c$ is a constant vector [5].

The following theorem and its result are the key in this article.

**Theorem 1.** Under the above notation, $\tilde{\gamma}$ is a Bertrand curve. Moreover, all Bertrand curves can be constructed by the above method [5].

**Corollary 1.** The spherical curve $\gamma$ is a circle if and only if the corresponding Bertrand curves are circular helices [5].

## 3. HELICES AND BERTRAND CURVES IN EUCLIDEAN 3-SPACE

In this section, we investigate Bertrand curves corresponding to the spherical images of the tangent indicatrix, binormal indicatrix, principal normal indicatrix and Darboux indicatrix of a space curve in Euclidean 3-space and give some results of helices and Bertrand curves. Moreover an example of the spherical image of the tangent indicatrix of a space curve and the corresponding Bertrand curve is presented.

**Theorem 2.** Let $\alpha$ be a space curve parametrized by arc-length parameter $s$. Then

$$\tilde{\gamma}_T(\sigma) = a\int_{\sigma_0}^{\sigma} T(\sigma)d\sigma + a\cot\theta \int_{\sigma_0}^{\sigma} B(\sigma)d\sigma + c \qquad (6)$$

is a Bertrand curve, where $a$, $\theta$ are constant numbers, $c$ is a constant vector and $\sigma = \int_0^s \|T'(s)\|ds.$

***Proof.*** We denote the spherical image of the tangent indicatrix of space curve $\alpha$ by $(T)$. By using Theorem 1, Bertrand curve corresponding to $(T)$ is

$$\tilde{\gamma}_T(\sigma) = a\int_{\sigma_0}^{\sigma} T(\sigma)d\sigma + a\cot\theta \int_{\sigma_0}^{\sigma} T(\sigma) \times T'(\sigma)d\sigma + c.$$

Applying the Frenet frame and the Frenet-Serret formulae, we obtain

$$\tilde{\gamma}_T(\sigma) = a\int_{\sigma_0}^{\sigma} T(\sigma)d\sigma + a\cot\theta \int_{\sigma_0}^{\sigma} B(\sigma)d\sigma + c.$$

This completes the proof.

Then, we have the following corollary of Theorem 2.

**Corollary 2.** Let $\alpha$ be a general helix parametrized by arc-length parameter $s$. Then

$$\tilde{\gamma}_T(\sigma) = a\int_{\sigma_0}^{\sigma} T(\sigma)d\sigma + a\cot\theta \int_{\sigma_0}^{\sigma} B(\sigma)d\sigma + c \tag{7}$$

is both a Bertrand curve and a circular helix.

***Proof.*** From Theorem 2, Bertrand curve corresponding to the spherical image of the tangent indicatrix $(T)$ of space curve $\alpha$ is

$$\tilde{\gamma}_T(\sigma) = a\int_{\sigma_0}^{\sigma} T(\sigma)d\sigma + a\cot\theta \int_{\sigma_0}^{\sigma} B(\sigma)d\sigma + c.$$

It is not difficult to show that $\alpha$ be a general helix if and only if the spherical image of the tangent $(T): I \to \mathbb{S}^2$ is a part of a circle on $\mathbb{S}^2$. Thus from Corollary 1, $\tilde{\gamma}_T$ is a circular helix.

**Theorem 3.** Let $\alpha$ be a space curve parametrized by arc-length parameter $s$. Then

$$\tilde{\gamma}_B(\sigma) = a\int_{\sigma_0}^{\sigma} B(\sigma)d\sigma + a\cot\theta \int_{\sigma_0}^{\sigma} T(\sigma)d\sigma + c \tag{8}$$

is a Bertrand curve, where $a$, $\theta$ are constant numbers, $c$ is a constant vector and $\sigma = \int_0^s \|B'(s)\|ds$.

***Proof.*** We denote the spherical image of the binormal indicatrix of space curve $\alpha$ by $(B)$. Hence applying Theorem 1, Bertrand curve corresponding to $(B)$ is

$$\tilde{\gamma}_B(\sigma) = a\int_{\sigma_0}^{\sigma} B(\sigma)d\sigma + a\cot\theta \int_{\sigma_0}^{\sigma} B(\sigma)\times B'(\sigma)d\sigma + c.$$

Using the Frenet frame and the Frenet-Serret formulae, we have

$$\tilde{\gamma}_B(\sigma) = a\int_{\sigma_0}^{\sigma} B(\sigma)d\sigma + a\cot\theta \int_{\sigma_0}^{\sigma} T(\sigma)d\sigma + c.$$

The proof is completed.
Thus we can give the following corollary of Theorem 3.

**Corollary 3.** Let $\alpha$ be a general helix parametrized by arc-length parameter $s$. Then

$$\tilde{\gamma}_B(\sigma) = a\int_{\sigma_0}^{\sigma} B(\sigma)d\sigma + a\cot\theta \int_{\sigma_0}^{\sigma} T(\sigma)d\sigma + c \tag{9}$$

is both a Bertrand curve and a circular helix.

***Proof.*** From Theorem 3, Bertrand curve corresponding to the spherical image of the binormal indicatrix $(B)$ of space curve $\alpha$ is

$$\tilde{\gamma}_B(\sigma) = a\int_{\sigma_0}^{\sigma} B(\sigma)d\sigma + a\cot\theta \int_{\sigma_0}^{\sigma} T(\sigma)d\sigma + c.$$

We can show that $\alpha$ be a general helix if and only if the spherical image of the binormal $(B): I \to \mathbb{S}^2$ is a part of a circle on $\mathbb{S}^2$. Then by Corollary 1, $\tilde{\gamma}_B$ is a circular helix.

**Theorem 4.** Let $\alpha$ be a space curve parametrized by arc-length parameter $s$. Then

$$\tilde{\gamma}_N(\sigma) = a\int_{\sigma_0}^{\sigma} N(\sigma)d\sigma + a\cot\theta \int_{\sigma_0}^{\sigma} C(\sigma)d\sigma + c \qquad (10)$$

is a Bertrand curve, where $a$, $\theta$ are constant numbers, $c$ is a constant vector and $\sigma = \int_0^s \|N'(s)\|ds.$

***Proof.*** We denote the spherical image of the principal normal indicatrix of space curve $\alpha$ by $(N)$. By using Theorem 1, Bertrand curve corresponding to $(N)$ is

$$\tilde{\gamma}_N(\sigma) = a\int_{\sigma_0}^{\sigma} N(\sigma)d\sigma + a\cot\theta \int_{\sigma_0}^{\sigma} N(\sigma) \times N'(\sigma)d\sigma + c.$$

Since $N(\sigma) \times N'(\sigma) = C(\sigma)$, we get

$$\tilde{\gamma}_N(\sigma) = a\int_{\sigma_0}^{\sigma} N(\sigma)d\sigma + a\cot\theta \int_{\sigma_0}^{\sigma} C(\sigma)d\sigma + c.$$

This completes the proof.

We can give the following results of Theorem 4.

**Corollary 4.** Let $\alpha$ be a slant helix parametrized by arc-length parameter $s$. Then

$$\tilde{\gamma}_N(\sigma) = a\int_{\sigma_0}^{\sigma} N(\sigma)d\sigma + a\cot\theta \int_{\sigma_0}^{\sigma} C(\sigma)d\sigma + c \qquad (11)$$

is both a Bertrand curve and a circular helix.

***Proof.*** From Theorem 4, Bertrand curve corresponding to the spherical image of the principal normal indicatrix $(N)$ of space curve $\alpha$ is

$$\tilde{\gamma}_N(\sigma) = a\int_{\sigma_0}^{\sigma} N(\sigma)d\sigma + a\cot\theta \int_{\sigma_0}^{\sigma} C(\sigma)d\sigma + c.$$

We can show that $\alpha$ be a slant helix if and only if the spherical image of the principal normal $(N): I \to \mathbb{S}^2$ is a part of a circle on $\mathbb{S}^2$. Then by Corollary 1, $\tilde{\gamma}_N$ is a circular helix.

**Corollary 5.** Let $\alpha$ be a general helix parametrized by arc-length parameter $s$. Then

$$\tilde{\gamma}_N(\sigma) = a\int_{\sigma_0}^{\sigma} N(\sigma)d\sigma + a\cot\theta C(\sigma) + c \qquad (12)$$

is a Bertrand curve.
*Proof.* By Theorem 4,

$$\tilde{\gamma}_N(\sigma) = a\int_{\sigma_0}^{\sigma} N(\sigma)d\sigma + a\cot\theta \int_{\sigma_0}^{\sigma} C(\sigma)d\sigma + c$$

is a Bertrand curve. Since $\alpha$ is a general helix, we can easily show that the Darboux indicatrix $C$ is constant. Thus we immediately have

$$\tilde{\gamma}_N(\sigma) = a\int_{\sigma_0}^{\sigma} N(\sigma)d\sigma + a\cot\theta C(\sigma) + c.$$

This completes the proof.

**Theorem 5.** Let $\alpha$ be a space curve parametrized by arc-length parameter $s$. Then

$$\tilde{\gamma}_C(\sigma) = a\int_{\sigma_0}^{\sigma} C(\sigma)d\sigma + a\cot\theta \int_{\sigma_0}^{\sigma} N(\sigma)d\sigma + c \qquad (13)$$

is a Bertrand curve, where $a$, $\theta$ are constant numbers, $c$ is a constant vector and $\sigma = \int_0^s \|C'(s)\|ds$.
*Proof.* We denote the Darboux indicatrix of space curve $\alpha$ by $(C)$. Then from Theorem 1, Bertrand curve corresponding to $(C)$ is

$$\tilde{\gamma}_C(\sigma) = a\int_{\sigma_0}^{\sigma} C(\sigma)d\sigma + a\cot\theta \int_{\sigma_0}^{\sigma} C(\sigma) \times C'(\sigma)d\sigma + c.$$

Since $C = \dfrac{\tau T + \kappa B}{\sqrt{\tau^2 + \kappa^2}}$, by taking the derivative of this equation with respect to $\sigma$, we have

$$C' = \frac{\kappa T - \tau B}{\sqrt{\tau^2 + \kappa^2}}. \qquad (14)$$

So we can easily obtain

$$C \times C' = N. \qquad (15)$$

Using the above equations, we can form

$$\tilde{\gamma}_C(\sigma) = a\int_{\sigma_0}^{\sigma} C(\sigma)d\sigma + a\cot\theta \int_{\sigma_0}^{\sigma} N(\sigma)d\sigma + c.$$

This completes the proof.

Now we give an example of the spherical image of the tangent indicatrix of a space curve and the corresponding Bertrand curve and draw their pictures by using Mathematica computer program.

**Example 1.** Let us consider a space curve $\alpha$ defined by

$$\alpha(s) = \left(-\cos s, \frac{\sin^2 s}{2}, \frac{1}{4}\sin 2s + \frac{s}{2}\right). \qquad (16)$$

The picture of the space curve $\alpha$ is given by Figure 1.

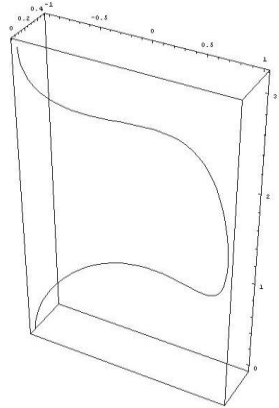

**Figure 1** - Space curve $\alpha$.

The spherical image of the tangent indicatrix of the space curve $\alpha$ is

$$\gamma(s) = (T) = (\sin s,\ \sin s \cos s,\ \cos^2 s). \qquad (17)$$

The picture of the spherical curve $\gamma$ is given by Figure 2.

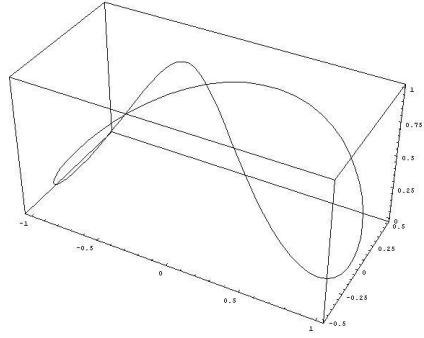

**Figure 2** - The spherical image of the tangent indicatrix $\gamma$ of the space curve $\alpha$.

From Theorem 2, we have the corresponding Bertrand curve ($a = 1$ and $\cot\theta = 1$):

$$\tilde{\gamma}_T(s) = \begin{pmatrix} -\cos s + \dfrac{1}{\sqrt{2}}\left( -\dfrac{1}{2}\int_0^s \dfrac{d\sigma}{\sqrt{1-\dfrac{1}{2}\sin^2\sigma}} + \dfrac{1}{2}\left( -4\int_0^s \sqrt{1-\dfrac{1}{2}\sin^2\sigma}\,d\sigma + 3\int_0^s \dfrac{d\sigma}{\sqrt{1-\dfrac{1}{2}\sin^2\sigma}} \right) \right), \\ -2\arctan\left( \dfrac{\sqrt{2}\sin s}{\sqrt{3+\cos(2s)}} \right) - \dfrac{1}{4}\cos(2s) + \dfrac{\sqrt{3+\cos(2s)}\sin s}{2\sqrt{2}}, \\ \dfrac{s}{2} - \dfrac{\cos s\sqrt{3+\cos(2s)}}{2\sqrt{2}} + \dfrac{3}{2}\log\left( \sqrt{2}\cos s + \sqrt{3+\cos(2s)} \right) + \dfrac{1}{4}\sin(2s) \end{pmatrix}$$

(18) [5].

We can draw the picture of the corresponding Bertrand curve $\tilde{\gamma}_T$ as follows (Figure 3).

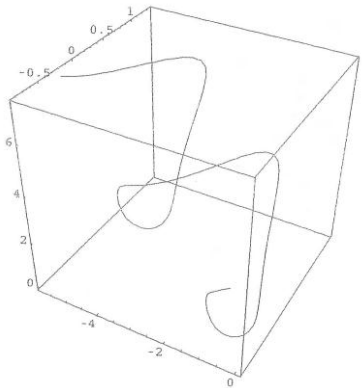

**Figure 3** - Bertrand curve $\tilde{\gamma}_T$ corresponding to the spherical image of the tangent indicatrix of the space curve $\alpha$ [5]